\documentclass[12pt]{article}
\usepackage{amsmath,amsfonts,amssymb,amsthm,epsfig}
\newcommand{\bb}{\mathbb}
\newcommand{\Q}{\bb Q}
\newcommand{\R}{\bb R}
\newcommand{\C}{\bb C}
\newcommand{\Z}{\bb Z}
\newcommand{\eps}{\varepsilon}
\newcommand{\hh}{\mathcal H}

\DeclareMathOperator{\area}{area}          
\DeclareMathOperator{\hol}{hol}
\DeclareMathOperator{\SL}{SL}

\newtheorem{theorem}{Theorem}[section]
\newtheorem{definition}[theorem]{Definition}
\newtheorem{corollary}[theorem]{Corollary}
\newtheorem{lemma}[theorem]{Lemma}
\newtheorem{proposition}[theorem]{Proposition}
\newtheorem{remark}[theorem]{Remark}

\begin{document}
\title{Minimal nonergodic directions on genus $2$ translation surfaces}
\author{
   Yitwah Cheung\\
   Northwestern University\\ Evanston, Illinois\\
   email: \emph{yitwah@math.northwestern.edu}
\bigskip\\ and \bigskip\\
   Howard Masur\thanks{This research is partially 
   supported by NSF grant DMS0244472.}\\
   University of Illinois at Chicago\\ Chicago, Illinois\\
   email: \emph{masur@math.uic.edu}
}
\maketitle

\begin{abstract}
It is well-known that on any Veech surface, the dynamics in 
any minimal direction is uniquely ergodic.  In this paper it 
is shown that for any genus $2$ translation surface which is 
not a Veech surface there are uncountable many minimal but 
not uniquely ergodic directions.  
\end{abstract}  

\section{Introduction}
Suppose $(X,\omega)$ is a translation surface where the genus of $X$ is 
at least $2$.  This means that $X$ is a Riemann surface, and $\omega$ is 
a holomorphic $1$-form on $X$.  For each $\theta\in[0,2\pi)$ there is a 
vector field defined on the complement of the zeroes of $\omega$ such 
that $\arg \omega=\theta$ along this vector field.   The corresponding 
flow lines are denoted $\phi_\theta$.  For a countable set of $\theta$ 
there is a flow line of $\phi_\theta$ joining a pair of zeroes of $\omega$.  
These flow lines are called saddle connections.  For any $\theta$ such 
that there is no saddle connection in direction $\theta$, the flow is 
minimal.  

Veech (\cite{V1}) found examples of certain skew rotations over 
a circle which are minimal, but not uniquely ergodic.  Namely, the 
orbits are dense but not uniformly distributed.  Veech's examples 
can be interpreted (\cite{MT}) in terms of flows on $(X,\omega)$ 
where $X$ has genus $2$ and $\omega$ has a pair of simple zeroes.  

Take two copies of the standard torus $\bb R^2/\bb Z^2$ and mark off 
a segment along the vertical axis from $(0,0)$ to $(0,\alpha)$, where 
$0<\alpha<1$.  Cut each torus along the segment and glue pairwise 
along the slits.  The resulting surface $(X_\alpha,\omega)$ is the
connected sum of the pair of tori and is a branched double cover over 
the standard torus, branched over $(0,0)$ and $(0,\alpha)$.  These two 
endpoints of the slits become the zeroes of order one of $\omega$.  
There are a pair of circles on $X_\alpha$ such that the first return 
map of $\phi_\theta$ to these circles gives a skew rotation over the 
circle.  If $\alpha$ is irrational, then there are directions $\theta$ 
such that the flow $\phi_\theta$ is minimal but not uniquely ergodic 
(\cite{MT}).  These reproduce the original Veech examples. 

Additional results about these examples are known.  Cheung (\cite{Ch1}) 
has shown that if $\alpha$ satisfies a Diophantine condition of the 
form that there exists $c>0,s>0$ such that $|\alpha-p/q|<c/q^s$ has 
no rational solutions $p/q$, then the Hausdorff dimension of the set 
of $\theta\in[0,2\pi)$ such that $\phi_\theta$ is not ergodic is 
exactly $1/2$.  On the other hand Boshernitzan showed (in an Appendix 
to the above paper) that there is a residual set of Liouville numbers 
$\alpha$ such that this set of $\theta$ has Hausdorff dimension $0$.  
The dimension $1/2$ in the Cheung result is sharp, for it was shown 
(\cite{M}) that for any $(X,\omega)$ (in any genus) this set of 
$\theta$ has Hausdorff dimension at most $1/2$.  

Now in the slit torus case, if $\alpha$ is rational, then minimality 
implies unique ergodicity.  This is part of a more general phenomenon 
called Veech dichotomy.  There is a natural action of $\SL_2(\R)$ on 
the moduli space of translation surfaces.  A \emph{Veech surface} is 
one whose stabilizer $\SL(X,\omega)$ is a lattice in $\SL_2(\R)$.  
These surfaces have the property that for any direction $\theta$, 
either the flow $\phi_\theta$ is periodic or it is minimal and 
uniquely ergodic (\cite{V2}).  

This raises the question of whether every surface $(X,\omega)$ that is 
not a Veech surface has a minimal but nonuniquely ergodic direction.  
In \cite{MS} it was shown that for every component of every moduli 
space of $(X,\omega)$ (other than a finite number of exceptional ones), 
there exists $\delta>0$ such that for almost every $(X,\omega)$ in 
that component, (with respect to the natural Lebesgue measure class) 
the Hausdorff dimension of the set of $\theta$ such that $\phi_\theta$ 
is minimal but not ergodic is $\delta$.  That theorem does not however 
answer the question for every surface. 

In this paper we establish the following converse to Veech dichotomy 
in genus $2$.  Let $\hh(1,1)$ be the moduli space of translation 
surfaces in genus $2$ with two simple zeroes and $\hh(2)$ the moduli 
space of translation surfaces with a single zero of order two.  
\begin{theorem}\label{thm:main}
For \emph{any} surface $(X,\omega)\in\hh(1,1)$ or $\hh(2)$ which is 
not a Veech surface,  there are an uncountable number of directions 
$\theta$ such that the flow $\phi_\theta$ is minimal but not uniquely 
ergodic.  
\end{theorem}
 
The theorem strengthens Theorem~1.5 of \cite{Mc2}.  There it was shown 
that for any surface in genus $2$ which is not a Veech surface, there 
is a direction $\theta$ such that the flow $\phi_\theta$ is not uniquely 
ergodic, and not all leaves are closed.  We also remark that in genus $2$, 
the Veech translation surfaces have been classified (\cite{Mc2}, \cite{Ca}).
We also remark that the converse to Veech dichotomy is false in higher genus. 
B.~Weiss (oral communication) using a construction of Hubert-Schmidt 
(\cite{zorich:Hubert:Schmidt:Infinitely:generated})
has provided an example in genus $5$ which is not a Veech surface and yet for which the Veech
dichotomy holds. The surface is a branched cover over a Veech surface in $\hh(2)$. 

\emph{Acknowledgments:} 
The authors are deeply indebted to the Yaroslav Vorobets for correcting 
an error in the original proof of Corollary~\ref{cor:irrational} and 
suggesting the more efficient proofs of Lemma~\ref{lem:3twists} and 
especially of Proposition~\ref{prop:key} which we have adopted over 
our original elementary, but significantly longer, arguments.

\section{Splittings in $\hh(2)$ and $\hh(1,1)$} 
We begin by generalising the construction of slit tori discussed 
in the introduction.  Suppose $T_1,T_2$ are a pair of flat tori 
defined by lattices $L_1,L_2$.  Let $l_1,l_2$ be simple segments 
on each, determining the same holonomy vector $w\in\C$.  Cut each 
$T_i$ along $l_i$ and glue the resulting tori together along the 
cuts.  The resulting surface is a connected sum of the pair of 
tori and belongs to the stratum $\hh(1,1)$.  The endpoints of the 
slits become simple zeroes.  If $(X,\omega)$ is biholomorphically 
equivalent to this surface, 
we say $(L_1,L_2,w)$ is a \emph{splitting} of $(X,\omega)$.  

Conversely, any surface in $\hh(1,1)$ can be constructed in this 
way for infinitely many possible triples $(L_1,L_2,w)$ (\cite{Mc1}).  

We can construct surfaces in $\hh(2)$ in a similar fashion.  
Given a lattice $L_1$ we may cut along a simple segment with 
holonomy $w$ as above, then identify opposite ends of the slit.  
This forms a torus with two boundary circles attached at a point.  
Glue in a cylinder, attaching a boundary component to each of the 
boundary circles.  The holonomy of the boundary circles is $w$.  
Every surface in $\hh(2)$ is found by such a construction (\cite{Mc1}).  
We again refer to $(L_1,L_2,w)$ as a splitting of $(X,\omega)$.  
In this case, we can think of the glued cylinder as a torus $T_2$ 
cut along a simple closed curve with holonomy $w$.  
We refer to this cylinder as a \emph{degenerate} torus.  

The result we use to construct nonuniquely ergodic directions 
is given in the following theorem (\cite{MS},\cite{MT}).  
\begin{theorem}\label{thm:nonuniq}
Let $(L_1^n,L_2^n,w^n)$ be a sequence of splittings of $(X,\omega)$ 
into tori $T_1^n,T_2^n$ and assume the directions of the vectors 
$w^n$ converge to some direction $\theta$.  Let $h_n>0$ be the 
component of $w^n$ in the direction perpendicular to $\theta$ 
and $a_n=\area(T_1^n\Delta T_1^{n+1})$ the area of the regions 
exchanged between consecutive splittings.  If 
\begin{enumerate}
  \item $\sum_{n=1}^\infty a_n<\infty$, 
  \item there exists $c>0$ such that $\area(T_1^n)>c$, $\area(T_2^n)>c$ 
        for all $n$, and 
  \item $\lim_{n\to\infty} h_n=0$, 
\end{enumerate}
then $\theta$ is a nonergodic direction.  
\end{theorem}

The idea behind the proof of Theorem~\ref{thm:main} is to construct 
uncountably many sequences of splittings satisfying the summability 
condition above with distinct limiting directions $\theta$.  Since 
there are only countably many nonminimal directions, there must be 
uncountably many limiting directions which are minimal but not ergodic.  

In Section~\ref{New} we find conditions to generate new 
splittings out of old and we find a useful estimate for 
the change in area of the corresponding tori.  Now if a 
surface $(X,\omega)$ is a Veech surface, then for any 
splitting $(L_1,L_2,w)$ of $(X,\omega)$, the vector $w$ 
is rational in each lattice $L_i$.  
In Section~\ref{Irrational} we show that if a surface is 
not Veech, there is some splitting $(L_1,L_2,w)$ such that 
the vector $w$ is irrational in one of the lattices $L_i$.  
We use this splitting to begin the inductive process of 
finding sequences of splittings that satisfy the hypotheses 
of Theorem~\ref{thm:nonuniq}.  The irrationality is used to 
make the change in area small and to ensure the inductive 
process can be continued.  The proof of this in the case of 
$\hh(2)$ is essentially the same although easier than in the 
case of $\hh(1,1)$.  We will focus on the latter and point 
out the differences with the former as they occur.

\section{Building New Splittings by Twists}\label{New}
We build new splittings out of old splittings by a Dehn twist 
operation.  The construction works in exactly the same way for 
the two moduli spaces.  

First let us adopt the notation that for vectors $v,v'\in\R^2$ 
given in rectangular coordinates by $v=\left<x_1,y_1\right>$ 
and $v'=\left<x_2,y_2\right>$ their cross product is 
$$v\times v':=x_1y_2-x_2y_1\in\R.$$  

Let $(L_1,L_2,w)$ be a splitting of $(X,\omega)$.  We may think 
of each slit torus $T_i$ as a closed subsurface (with boundary) 
in $X$ separated by a pair of saddle connections $\alpha_{\pm}$ 
given by $$\partial T_1 = \alpha_+ -\alpha_-,\quad
                              \hol(\alpha_+)=w=\hol(\alpha_-).$$  
We would like to have an explicit construction of a class of 
simple closed curves having nonzero (geometric) intersection 
with $\alpha_+\cup\alpha_-$.  Let $C_i$ be a cylinder contained 
in $T_i$ that is disjoint from the line segment $\ell_i$.  
(In the case where $T_2$ is degenerate we allow $C_2$ to consist 
of a single closed curve.)  Furthermore, assume that the holonomy 
$v_i$ of the core curve $\gamma_i$ is not parallel to $w$.  
Each $\gamma_i$ has a unique translate crossing the midpoint 
of $\ell_i$.  As closed curves joining the boundary of $T_i$ 
in $X$ to itself, these translates can be concatenated to form 
a simple closed curve $\gamma$ whose holonomy is $v_1+v_2$.  
Note that for $\gamma$ to be well-defined the curves $\gamma_1$ 
and $\gamma_2$ must have compatible orientations, i.e. 
\begin{equation}\label{eq:same:side}
               (v_1\times w)(v_2\times w)>0.  
\end{equation}
Also, since $T_i\setminus C_i$ is a parallelogram with sides 
given by the vectors $v_i$ and $w$, we have 
\begin{equation}\label{eq:disjoint}
     |v_i\times w|\le A_i=\area(T_i)\quad\text{for $i=1,2$.}  
\end{equation}
Conversely, suppose $v_i\in L_i,i=1,2$ are primitive vectors 
such that (\ref{eq:same:side}) and (\ref{eq:disjoint}) hold.  
The latter condition guarantees that $v_i$ is the holonomy of 
a simple closed curve $\gamma_i$ on $T_i$ which can be realized 
by a saddle connection joining the initial point of the slit to 
itself, and another (provided $T_i$ is not degenerate) joining 
the terminal point to itself.  Neither of these intersect the 
interior of the slit.  The pair of saddle connections bound a 
closed cylinder $C_i$ in $T_i$ of curves which do not cross the 
slit, where we allow $C_2$ to be a single closed curve if $T_2$ 
is degenerate.  The complement of (the union of) these cylinders 
in $X$ is an open annulus which, by (\ref{eq:same:side}), has a 
core curve $\gamma$ with holonomy $v_1+v_2$.  Thus, any pair of 
primitive vectors $(v_1,v_2)\in L_1\times L_2$ that satisfies 
(\ref{eq:same:side}) and (\ref{eq:disjoint}) uniquely determines 
a unique simple closed curve $\gamma=\gamma(v_1,v_2)$.  

For the purpose of geometric intuition, we shall give a normal 
form representation of the surface $(X,\omega)$ as an even-sided 
polygon in the plane with pairs of parallel sides identified.  
Each $T_i$ decomposes into $C_i$ and a parallelogram $R_i$, 
whose sides are given by the vectors $v_i$ and $w$.  Using the 
action of $\SL_2(\R)$ we may represent $R_1$ as a rectangle in 
the plane with sides parallel to the coordinate axes.  The 
parallelogram $R_2$ may be placed adjacent to the right edge 
of $R_1$ (i.e. the slit $\alpha_+$), while a suitable choice 
of parallelograms representing the cylinders $C_1$ and $C_2$ 
may be placed adjacent to the lower edges of $R_1$ and $R_2$, 
respectively, so that there is no overlap.  
\begin{figure}[h]%%%%%%%%%%%%%%%%%%%%%%%%%%%%%%%%%%%%%%
%%%%%%%%%%%%%%%
\begin{center}
\includegraphics{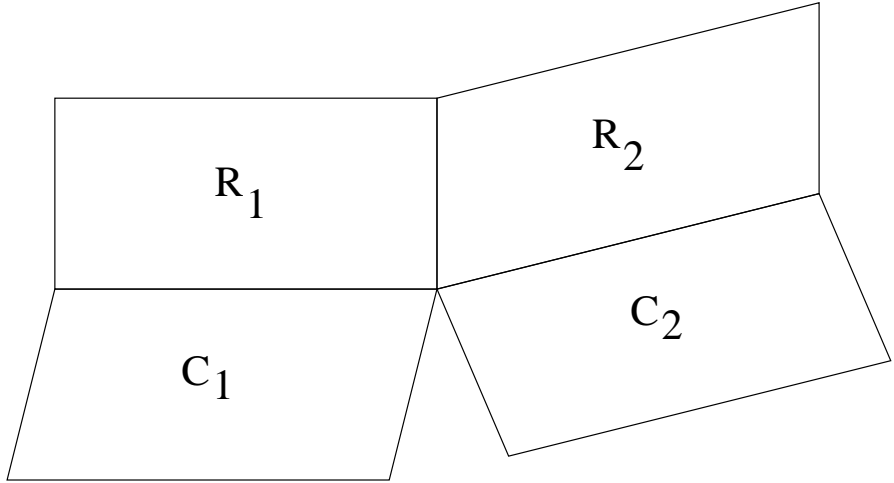}
\caption{Normal form for $(X,\omega)$.}\label{fig:normal}  
\end{center}
\end{figure}%%%%%%%%%%%%%%%%%%%%%%%%%%%%%%%%%%%%%%%%
%%%%%%%%%%%%%%%%%%
The boundary identifications are as follows.  The right edge of 
$R_2$ is identified with the left edge of $R_1$ to form the slit 
$\alpha_-$ homologous to $\alpha_+$.  For each $i=1,2$, the upper 
edge of $R_i$ is identified with the lower edge of $C_i$, while 
the left edge of $C_i$ is identified with its right edge.  In the 
case that $(X,\omega)$ has a double zero the cylinder forming the 
degenerate torus $T_2$ is represented as a parallelogram $R_2$ with 
upper and lower edges identified.  (Recall that in this case $C_2$ 
is a closed curve and not a cylinder.)  In both cases, the open 
annulus that forms the interior of $R_1\cup R_2$ has $\gamma$ 
as its core curve.  

Now consider the curves $\beta_{\pm}^k$ obtained by twisting 
the slits $\alpha_{\pm}$ (relative to their endpoints) $k$ times 
about $\gamma$ in the positive sense, i.e. right-twist if $k>0$, 
and left-twist if $k<0$.  In general, the geodesic representative 
of a twisted curve is a finite sequence of saddle connections.  
Let $w^k$ be the common holonomy vector of $\beta_{\pm}^k$.  Note 
that if $v_1\times w>0$, i.e. $\gamma$ is \emph{positively oriented} 
with respect to $\alpha_\pm$, then $w^k=w+k(v_1+v_2)$.  
\begin{lemma}\label{lem:twist}
The twisted curves $\beta_+^k$ and $\beta_-^k$ are simultaneously 
realized by a (single) saddle connection if and only if $w$ and 
$w^k$ lie on the same side of $v_1$ and $v_2$.  Thus, if $\gamma$ 
is positively oriented with respect to $\alpha_{\pm}$, then the 
condition is $$v_1\times w^k>0\quad\text{and}\quad v_2\times w^k>0$$ 
while both inequalities are reversed in the case where $\gamma$ is 
negatively oriented with respect to $\alpha_{\pm}$.  
\end{lemma}
\begin{proof}
First, recall that the Dehn twist operation depends only on the 
orientation of the surface and not that of the curves $\alpha_{\pm}$ 
and $\gamma$.  In particular, we may assume the orientation of 
$\gamma$ is chosen positive with respect to $\alpha_\pm$.  Since 
the hyperelliptic involution interchanges $\alpha_+$ and $\alpha_-$ 
while fixing $\gamma$, 
$\beta_+^k$ is realised by a saddle connection if and only if 
$\beta_-^k$ is.  If  $v_1\times v_2\ge0$ we show   $\beta_-^k$ 
is realized.  Otherwise we consider  $\beta_+^k$.  Without loss of generality then assume
$v_1\times v_2\geq 0$.   Since the action of 
$\SL_2(\R)$ preserves cross products, it is enough to verify that 
for \emph{some} $g\in\SL_2(\R)$, the vector $gw^k$ is the holonomy 
of a curve realised by a saddle connection on the flat surface 
$g\cdot(X,\omega)$.  Consider first the case of a left twist $k<0$.  
Choose $g$ so that $v_1$ is horizontal and $w$ is vertical, 
pointing upwards as in Figure~\ref{fig:twist}.  
\begin{figure}[h]%%%%%%%%%%%%%%%%%%%%%%%%%%%%%%%%%%%%%%
%%%%%%%%%%%%%%%
\begin{center}
\includegraphics{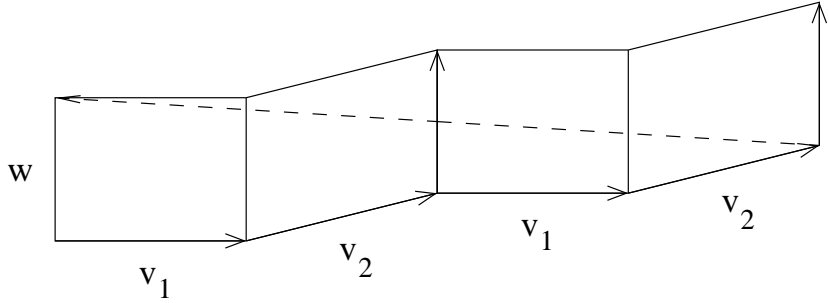}
\caption{$\beta_-^{-2}$ as a single saddle connection.}  
\label{fig:twist}
\end{center}
\end{figure}%%%%%%%%%%%%%%%%%%%%%%%%%%%%%%%%%%%%%%%%
%%%%%%%%%%%%%%%%%%
One sees easily that $\beta_-^k$ is realised by a saddle connection 
if and only if $-k$ times the vertical component of $v_2$ is less 
than $|w|$, which is equivalent to $v_1\times w^k>0$. 
 On the other 
hand, $k<0$ implies $v_2\times w^k=v_2\times w + kv_2\times v_1>0.$

The case of a right twist is similar.  The condition $k>0$ implies 
$v_1\times w^k>0$.  By choosing $g$ so that $v_2$ is horizontal and 
$w$ is vertical we see that $v_2\times w^k>0$ is equivalent to the 
condition that $\beta_-^k$ is realised by a saddle connection.  
\end{proof}

\begin{definition}
Given a splitting $(L_1,L_2,w)$ and a pair of primitive vectors 
$(v_1,v_2)\in L_1\times L_2$ satisfying (\ref{eq:same:side}) and 
(\ref{eq:disjoint}) we say $v_2$ is a \emph{good partner} for $v_1$ 
with $n$ twists if 
\begin{equation}\label{eq:n:twists}
  |v_1\times v_2|<\frac{1}{n}\max(|v_1\times w|,|v_2\times w|).  
\end{equation}
Note that if $v_2$ is a good partner for $v_1$ with $n$ twists, 
then Lemma~\ref{lem:twist} implies there are $n$ splittings 
$(L_l^k,L_2^k,w^k)$ of $(X,\omega)$ with $w^k=w+k(v_1+v_2)$, 
where $k=1,\dots,n$ or $k=-1,\dots,-n$.  
\end{definition}

Next, we estimate the total area of the regions exchanged between 
the new and old splittings. 
\begin{lemma}
\label{lem:area}
 Let $(L_1',L_2',w')$ denote the new 
splitting of $(X,\omega)$ into tori $T_1',T_2'$, determined by a 
particular value of $k, |k|\le n$. 
Then $T_1\Delta T_1'=T_2\Delta T_2'$ and 
\begin{equation}
\label{eq:area}
  \area(T_1\Delta T_1') \le  |v_1\times w|+|v_1\times w'|
                        \le 2|v_1\times w|+n|v_1\times v_2|.  
\end{equation}
\end{lemma}
\begin{proof}
The first statement is clear.
For the second we note that the torus $T_1'$ contains $C_1$ 
and $T_2'$ contains $C_2$.   Since $C_1\subset T_1\cap T_1'$, it 
follows that 
$T_1\Delta T_1'\subset(T_1\setminus C_1)\cup(T_1'\setminus C_1)$.  
By construction, $\area(T_1\setminus C_1)=|v_1\times w|$ and 
$\area(T_1'\setminus C_1)=|v_1\times w'|$. 
 This gives the desired inequality.  
\end{proof}

\section{Irrational Splittings}\label{Irrational}
To construct new splittings with small exchange of area, 
we need the notion of an irrational splitting.  
\begin{definition}
We say the splitting $(L_1,L_2,w)$ is \emph{rational} in $L_i$ 
if the vector $w$ is a scalar multiple of some element in $L_i$; 
otherwise, it is \emph{irrational} in $L_i$.  The splitting is 
\emph{irrational} if it is irrational in either $L_1$ or $L_2$.  
\end{definition}
Note that if $(X,\omega)$ has a splitting that is rational in 
$L_i$  then the flow in the direction of $w$ is 
 periodic in the torus $T_i$. 

In Proposition~\ref{prop:key} we will assume that the splitting 
is irrational in at least one of the two lattices.  (In the case 
of $(X,\omega)$ having a double zero, we will assume that the 
splitting is irrational in $L_1$.)  We will then find a new 
splitting so that the right hand side of (\ref{eq:area}) is small.  
To continue the process of finding additional splittings we will 
need to know that the new splitting is still irrational.  That 
fact will be accomplished by the next lemma.  
\begin{lemma}\label{lem:3twists}
Let $(L_1,L_2,w)$ be a splitting that is irrational in $L_1$ and 
suppose $v_1\in L_1$ has a good partner $v_2\in L_2$ with three 
twists.  Then at least one of the three twists gives a splitting 
$(L_1',L_2',w')$ that is irrational in $L_1'$.  
\end{lemma}
\begin{proof}
Let $\gamma_0$ be a simple segment in $C_1$ that concatenates 
with $\alpha_-$ to form a simple closed curve in $T_1$.  Then 
the lattice $L_1$ is generated by $v_0+w$ and $v_1$ where $v_0$ 
is the holonomy of $\gamma_0$.  Observe that $(L_1,L_2,w)$ is 
rational in $L_1$ if and only if $w$ is a scalar multiple of 
a vector in the lattice $L_0$ generated by $v_0$ and $v_1$.  
Indeed, $w = c(av_0+bv_1)$ for some $a,b\in\Z$ and $c>0$ if 
and only if $(1+ac)w = c(a(v_0+w)+bv_1)$ where $1+ac\neq0$ 
 since $v_0+w$ and $v_1$ are linearly independent.  Now since 
$T_1$ and $T_1'$ share the same cylinder $C_1$, it follows 
that the splitting $(L_1',L_2',w')$ is rational in $L_1'$ 
if and only if $w'\in L_0$.  Thus, if $(L_1',L_2',w')$ with 
$w'=w+k(v_1+v_2)$ is rational in $L_1'$ for $k=1,2,3$ or 
$k=-1,-2,-3$, then the three vectors $w'$ are parallel to 
elements in $L_0$.  Let $L_0'$ be the lattice generated by 
$w$ and $v_1+v_2$.  Since $L_0$ and $L_0'$ are lattices in 
$\R^2$ whose vectors share three nonparallel directions, 
they are isogenous, (cf.
proof of Theorem 7.3 of \cite{Mc2})  i.e. $L_0'\subset \lambda L_0$ for some 
$\lambda>0$.  This implies that in fact they share all possible directions and so the splitting
$(L_1,L_2,w)$ is 
rational in $L_1$.  
\end{proof}

The next proposition will help us find some initial splitting 
which is irrational so that the inductive process can begin.  
\begin{proposition}\label{prop:new}
Let $(L_1,L_2,w)$ be a splitting of $(X,\omega)$, rational in 
both $L_1$ and $L_2$.  Suppose that there is a $v_1\in L_1$ 
which has two good partners $v_2,v_2'\in L_2$ with at least $2$ 
twists each and such that $v_2-v_2'$ is not parallel to $w$.  
Then, if $(X,\omega)$ is a not Veech surface, at least one 
of the four twists gives a splitting $(L_1',L_2',w')$ that 
is irrational in $L_1'$.  
\end{proposition}
\begin{proof}
Since the property of irrationality is invariant under the 
$\SL_2(\bb R)$ action and scaling, we may normalize so that 
$v_1=(1,0)$ and $w=(0,1)$.  In terms of the normal form, the 
parallelogram $R_1$ is a unit square with lower right corner 
at the origin.  Let $v_0$ be the holonomy of $\gamma_0$ as in 
the proof of the previous lemma.  The rationality of $(L_1,L_2,w)$ 
in $L_1$ and $L_2$ is equivalent to 
\begin{equation}\label{eq:rat:twists}
   v_0\times w\in\Q, \quad\text\quad \frac{v_2\times w}{v_2'\times w}\in\Q.  
\end{equation}
Let $\theta:=v_0\times v_1$.  Note that $-\theta$ is the 
height of the cylinder $C_1$; in particular, $\theta\neq0$.  
The slope of the vector $w^k=w+k(v_1+v_2)$ is given by 
\begin{equation}\label{eq:slope}
   \sigma_k = \frac{-v_1\times w^k}{w\times w^k} 
            = \frac{1+k(v_1\times v_2)}{k(1+v_2\times w)}.  
\end{equation}
The first return map to the base of $R_1$ under the flow in 
direction $w^k$ is a rotation with rotation number $\rho_k$.  
Note that a point starting at the base of $R_1$ returns to 
a point at the top of $R_1$ with the same $x$-coordinate.  
The identification of the top of $R_1$ with the base of $C_1$ 
shifts the $x$-coordinate by $-v_0\times w$.  Provided that 
$\sigma_k\neq0$, the total shift in the $x$-coordinate for 
the first return to the base of $R_1$ is given by 
\begin{equation}\label{eq:rotation}
   \rho_k = -\frac{\theta}{\sigma_k}-v_0\times w.  
\end{equation}
Since $v_2$ is a good partner of $v_1$ with two twists, 
we have $\rho_k\in\Q$ for $\pm k=1,2$.  Thus, together 
with the first part of (\ref{eq:rat:twists}) we have 
\begin{equation}\label{eq:return}
   \sigma_k\neq0 \quad\Rightarrow\quad 
   \frac{\theta(1+v_2\times w)}{1+k(v_1\times v_2)}\in\Q.  
\end{equation}
If $v_1\times v_2\neq0$ and $\sigma_k\neq0$ for both values 
of $k$ then taking ratios in (\ref{eq:return}) we see that 
\begin{equation}\label{eq:rat:prods}
   v_1\times v_2\in\Q,\quad\text\quad\theta(1+v_2\times w)\in\Q.  
\end{equation}
On the other hand, if either $v_1\times v_2=0$ or $\sigma_k=0$ 
then $v_1\times v_2\in\Q$, and using (\ref{eq:return}) for the 
other value of $k$ we see that (\ref{eq:rat:prods}) still holds.  
The same argument can be applied to $v_2'$ to yield 
\begin{equation}\label{eq:rat:prods'}
   v_1\times v_2'\in\Q,\quad\text\quad\theta(1+v_2'\times w)\in\Q 
\end{equation}
so that taking ratios again, we have 
\begin{equation}\label{eq:rat:ratio}
   \frac{1+v_2\times w}{1+v_2'\times w}\in\Q.  
\end{equation}
Since $v_2'-v_2$ is not parallel to $w$, the above together 
with the second part of (\ref{eq:rat:twists}) implies both 
\begin{equation}\label{eq:rat:x's}
   v_2\times w, \quad\text\quad v_2'\times w\in\Q.  
\end{equation}
Hence, $\theta\in\Q$, which together with (\ref{eq:rat:x's}) and 
the first parts of (\ref{eq:rat:twists}), (\ref{eq:rat:prods}), and (\ref{eq:rat:prods'}) 
implies $(X,\omega)$ is a branched cover of  the standard torus, branched over rational points. 
By \cite{GJ}, such a 
surface is Veech.
\end{proof}

The property that a surface $(X,\omega)$ admits a splitting 
$(L_1,L_2,w)$ together with vectors $v_1,v_2,v_2'$ satisfying 
the requirements in Proposition~\ref{prop:new} defines an 
open subset in $\hh(1,1)$.  (It is not hard to show that this 
property does not hold for \emph{any} surface in $\hh(2)$.)  
We show that this set is nonempty by constructing an explicit 
example: Let $v_1=(1,0)$, $w=(0,1)$, $v_2=(3,-1)$, $v_2'=(4,1)$, 
$L_1=\Z\times2\Z$ and $L_2$ the lattice generated by $v_2$ and 
$v_2'$.  Then $A_2=7$ and the fact that both $v_2$ and $v_2'$ 
have slopes strictly between $\pm1/2$ implies each has $2$ twists.  

Using results of McMullen we now obtain the following.  
\begin{corollary}\label{cor:irrational}
Suppose  $(X,\omega)$ is not a Veech surface. 
Then it admits an irrational splitting.  
\end{corollary}
\begin{proof}
Case 1: $(X,\omega)$ is an eigenform for real multiplication.  
(We refer the reader to \cite{Mc1} for the definition 
of an \emph{eigenform form real multiplication}.)  
The Corollary is given by Theorem~7.5 of \cite{Mc2}. 
 
Case 2: $(X,\omega)\in\hh(2)$.  Note that every cylinder 
determines a splitting.  If every splitting is rational, 
then $(X,\omega)$ has the property that every cylinder 
belongs to a cylinder decomposition; by Theorem 7.3 of 
\cite{Mc2} this means $(X,\omega)$ is an eigenform for 
real multiplication and in $\hh(2)$ this implies that the surface is a Veech surface.

Case 3: $(X,\omega)\in\hh(1,1)$ is not an eigenform for 
real multiplication.  By \cite{Mc1} the $\SL_2(\R)$-orbit 
comes arbitrarily close (up to scale) to the example 
described after Proposition~\ref{prop:new}, which now 
implies $(X,\omega)$ admits an irrational splitting.  
\end{proof}

\begin{remark}
Corollary~\ref{cor:irrational} shows that the 
hypothesis in Theorem~7.5 of \cite{Mc2} is unnecessary.  
\end{remark}

Theorem~\ref{thm:main} is essentially a consequence 
of the preceding Corollary and the next Proposition, 
followed by an application of Theorem~\ref{thm:nonuniq}.  
\begin{proposition}\label{prop:key}
If $(L_1,L_2,w)$ be an irrational splitting of $(X,\omega)$ 
into tori $T_1,T_2$, then for any $\eps>0$ there exists a 
\emph{new} irrational splitting $(L_1',L_2',w')$ into tori 
$T_1',T_2'$ such that the angle between $w,w'$ is 
$<\eps$ and $\area(T_1\Delta T_1')<\eps$.  
\end{proposition}

To prove Proposition~\ref{prop:key} we adopt the suggestion 
of Yaroslav Vorobets and exploit a theorem of McMullen whose 
proof itself uses Ratner's theorem.  This proof replaces our 
original more elementary but significantly longer proof. 

Let $G=\SL_2(\R)$ and $\Gamma=\SL_2(\Z)$ and regard $G/\Gamma$ as 
the space of oriented lattices $\Lambda\subset\R^2$ of coarea $1$.  
Let $N\subset G$ be the subgroup preserving horizontal vectors.  
Let $G_\Delta=\{(g,g)|g\in G\}$ and $N_\Delta=\{(g,g)|g\in N\}$.  
\begin{theorem}{\bf [Mc1,Theorem~2.6]}\label{thm:ratner}
Let $z=(\Lambda_1,\Lambda_2)\in(G\times G)/(\Gamma\times\Gamma)$ 
be a pair of lattices, and let $Z=\overline{N_\Delta z}$.  Then exactly 
one of the following holds.  
\begin{enumerate}
  \item There are horizontal vectors $v_i\in\Lambda_i$ with 
        $|v_1|/|v_2|\in\Q$.  Then $Z=N_\Delta z$.  
  \item There are horizontal vectors $v_i\in\Lambda_i$ with 
        $|v_1|/|v_2|\not\in\Q$.  Then $Z=(N\times N)z$.  
  \item One lattice, say $\Lambda_2$, contains a horizontal 
        vector but the other does not.  Then $Z=(G\times N)z$.  
  \item Neither lattice contains a horizontal vector, but 
        $\Lambda_1\cap\Lambda_2$ is of finite index in both.  
        Then $Z=G_\Delta z$.  
  \item The lattices $\Lambda_1$ and $\Lambda_2$ are 
        incommensurable, and neither contains a horizontal vector.  
        Then $Z=(G\times G)z$.  
\end{enumerate}
\end{theorem}
\begin{proof}[Proof of Proposition~\ref{prop:key}]
Let $(L_1,L_2,w)$ be a splitting of $(X,\omega)$ that is 
irrational in $L_1$.  Let $A_i$ be the coarea of $L_i$.  
It is enough to show that for any $\eps>0$ and sufficiently 
small $\eps'<\eps$ there is a pair 
of primitive vectors $(v_1,v_2)\in L_1\times L_2$ satisfying 
\begin{enumerate}%inserted (i) and (ii) labels
  \item[(i)] $|v_1\times w|<\min(\eps'/4,A_1)$ and $|v_2\times w|\le A_2$, 
  \item[(ii)] $|v_1\times v_2|<\eps'/6$ and 
        $|v_1\times v_2|<\frac{1}{3}\max(|v_1\times w|,|v_2\times w|)$.  
\end{enumerate}
Indeed, one of the three splittings obtained by twisting about 
$\gamma(v_1,v_2)$ will be irrational by Lemma~\ref{lem:3twists}. 
By Lemma~\ref{lem:area} $$\area(T_1\Delta T_1')\leq 2\eps'/4+3\eps'/6<\eps$$
and by choosing a 
sufficiently small $\eps'$, the direction of the vectors $v_1,v_2$, 
and hence that of $w+k(v_1+v_2)$, can be made within $\eps$ of $w$.  

% (This sentence is not necessary.)We shall regard $(G\times G)/(\Gamma\times\Gamma)$ %
%asthe space of pairs $(\Lambda_1,\Lambda_2)$ of oriented lattices in $\R^2$% 
%such that the coarea of $\Lambda_i$ is $A_i$ for $i=1,2$.%  
Let $Z$ be the closure of the orbit of $(L_1,L_2)$ under the 
action of $N_\Delta$ where $N\subset G$ is the subgroup that 
preserves vectors parallel to $w$.

Suppose first that $(L_1,L_2,w)$ is irrational in both $L_1$ and $L_2$.  
Then we are in case 4 or  5 of Theorem~\ref{thm:ratner}.  If case 5 holds, 
then the $G\times G$ action allows one to conclude that for any 
$\eps'>0$ there exist $(\Lambda_1,\Lambda_2)\in Z$ and primitive 
vectors $v_i'\in\Lambda_i$ both perpendicular to $w$ and having 
lengths $$|v_i'|<\frac{1}{|w|}\min(\eps',A_i).$$  If case 4 holds, then $\Lambda_1$ and
$\Lambda_2$ are isogenous and so once one %inserted "one" 
has found a vector $v_1'\in \Lambda_1$ perpendicular
to $w$ by the $G$ action on lattices,  then automatically there is a vector $v_2'\in \Lambda_2$
perpendicular to $w$ as well. Then we can use the $G_\Delta$ action to make both vectors
satisfy the above inequality as well.  

It now follows that there exists $g_n\in N$ and $v_{n,i}\in L_i$ such that $g_n v_{n,i}\to
v_i'$ as $n\to\infty$. 
Since $g\in G$ preserves cross products, and $g\in N$ fixes $w$, for large $n$ we have 
vectors $v_{n,1},v_{n,2}$ that satisfy
\begin{itemize}
\item 
$|v_{n,1}\times w|=|g_n v_{n,1}\times g_n w|=|g_n v_{n,1}\times w|<\min (\eps'/4,A_1)$
\item $|v_{n,2}\times w|=|g_n v_{n,2}\times g_n w|=|g_n v_{n,2}\times w|<A_2$
\item 
$|v_{n,1}\times v_{n,2}|=|g_n v_{n,1}\times g_n v_{n,2}|<\eps'/6$
\item $|v_{n,1}\times v_{n,2}|= |g_n v_{n,1}\times g_n v_{n,2}|<\frac{1}{3} \max(|g_n
v_{n,1}\times w|,|g_n v_{n,2} \times w|)=\frac{1}{3}\max(|v_{n,1}\times w|,|
v_{n,2}\times
w|)$.
\end{itemize}
This means that the vectors $v_{1,n},v_{2,n}$ satisfy the desired conditions (i) and (ii). 

Suppose next that $(L_1,L_2,w)$ is rational in $L_2$.  We are assuming 
that the splitting $(L_1,L_2,w)$ is irrational in $L_1$.  Now case 3 
of Theorem~\ref{thm:ratner} implies for any $\eps'>0$ there exist 
$(\Lambda_1,\Lambda_2)\in Z$ 
and primitive vectors $v_i'\in\Lambda_i$ both perpendicular to $w$ 
with $|v_1'|<\frac{1}{|w|}\min(\eps',A_1)$ and $|v_2'\times w|=A_2$.  
Again, we approximate $(\Lambda_1,\Lambda_2)$ by pairs of lattices 
in the $N$-orbit of $(L_1,L_2)$ and note that $|v_2\times w|= A_2$
for any pair $(v_1,v_2)$ that approximates $(v_1',v_2')$.  
\end{proof}

\begin{proof}[Proof of Theorem~\ref{thm:main}]
Suppose $(X,\omega)$ is not a Veech surface.  Since there 
are only countably many nonminimal directions, it suffices 
to construct an uncountable number of nonergodic directions.  
By Corollary~\ref{cor:irrational}, there is some irrational 
splitting $(L_1,L_2,w)$ of $(X,\omega)$ into tori $T_1,T_2$.  
Let $A_i$ be the area of $T_i$.  We inductively construct an 
infinite binary tree of splittings $(L_1^n(j),L_2^n(j),w^n(j))$ 
of $(X,\omega)$ so that there are $2^n$ splittings at  level 
$n$, which we index by $j=1,\dots,2^n$. We take the 
unique splitting of level zero to be $(L_1,L_2,w)$.  
At the completion of the $n$th level, we define $\eps_n>0$ 
to be the minimum distance in angle between the directions 
of any two splittings that have been constructed up to this 
point.  To construct the splittings of the next level, we 
apply Proposition~\ref{prop:key} with $\eps<\eps_n/4$.   
Continuing ad infinitum we obtain an uncountable number of 
sequences $(L_1^n,L_2^n,w^n)$ of splittings of $(X,\omega)$.
Since $\eps_n\leq (1/2)^n$ the directions of the $\{w^n\}$ for any infinite geodesic in the tree
converge to
a limiting direction $\theta$.    Let $h_n$ be as in Theorem~\ref{thm:nonuniq}.  
Then   
$$h_{n+1} = |w^{n+1}|\sin(\angle w^{n+1}\theta)
        \le |w^{n+1}|\sin(2\angle w^nw^{n+1}) 
        \le \frac{2|w^n\times w^{n+1}|}{|w^n|}.$$ 
We have $\lim|w^n|=\infty$.  Now    
$|w^n\times(v_1+v_2)|\leq A_1+A_2$ since the left side is the area of an  annulus 
contained in $(X,\omega)$. Since   
 $w^{n+1}=w^n+k(v_1+v_2)$ for some $k,|k|\leq 3,$
$$|w^n\times w^{n+1}|\le3(A_1+A_2).$$
As a result we  have $\lim_{n\to\infty} h_n=0$.  The minimum spacings $\eps_n$ at level
$n$ satisfy $$\eps_{n+1}\leq \eps_n/2.$$
 Since the change in areas is at most $\eps_n$, the sum of the change of areas is at most  
$2\eps_0$.   If we choose $\eps_0<\frac{1}{2}\min(A_1,A_2)$ the remaining 
hypothesis of 
Theorem~\ref{thm:nonuniq} is 
satisfied.  

That the  nonergodic limiting directions that have 
been constructed are all distinct follows from the condition that the spacing between a splitting at
level $n$ and either of its descendents is at most $\eps_n/4$. 
This completes the proof of Theorem~\ref{thm:main}.  
\end{proof}


\begin{thebibliography}{999}

\bibitem[Ca]{Ca} 
  K.~Calta, 
  {\em Veech surfaces and complete periodicity in genus two}, 
      to appear, JAMS.  

\bibitem[Ch1]{Ch1} 
  Y.~Cheung, 
  {\em Hausdorff dimension of the set of nonergodic directions}, 
      Ann. of Math., {\bf 158} (2003), 661-678.  

\bibitem[Ch2]{Ch2}
  Y.~Cheung,
  {\em Slowly divergent trajectories in moduli space}, 
      Conform. Geom. Dyn., {\bf 8} (2004), 167-189.  

\bibitem[GJ]{GJ}
  E.~Gutkin and C.~Judge,
  {\em Affine mappings of translation surfaces: geometry and arithmetic},
      Duke Math. J. {\bf 103} (2000), 191-213.  

\bibitem[HS]{zorich:Hubert:Schmidt:Infinitely:generated}
  P.~Hubert and T.~A.~Schmidt,
  {\em Infinitely generated Veech groups},
      Duke Math. J. {\bf 123} (2004), 49--69.

\bibitem[Mc1]{Mc1} 
  C.~McMullen, 
  {\em Dynamics of $\SL_2(\bb R)$ over moduli space in genus two},
     preprint, April 26, 2004.  
   
\bibitem[Mc2]{Mc2} 
  C.~McMullen,
  {\em Teichm\"uller curves in genus two: The decagon and beyond},
     preprint, May 11, 2004.  
   
\bibitem[M]{M} 
  H.~Masur, 
  {\em Hausdorff dimension of the set of nonergodic foliations of a
     quadratic differential} Duke Math. Journal {\bf 66} (1992) 387-442.  

\bibitem[MS]{MS}
  H.~Masur and J.~Smillie,
  {\em Hausdorff dimension of sets of nonergodic measured foliations},
     Ann. of Math., {\bf 134} (1991), 455-543.

\bibitem[MT]{MT} 
  H.~Masur and S.~Tabachnikov, 
  {\em Rational Billiards and Flat Structures},
     Handbook of Dynamical Systems, B.~Hasselblatt, A.~Katok eds. 
     Elsevier (2002), 1015-1089.  

\bibitem[V1]{V1} 
  W.~Veech, 
  {\em Strict ergodicity in zero dimensional dynamical systems 
      and the Kronecker-Weyl theorem mod 2}, 
      Transactions Amer. Math Soc. {\bf 140} (1969) 1-34.  
  
\bibitem[V2]{V2} 
  W.~Veech, 
  {\em Teichmuller curves in moduli space, Eisenstein series and 
      an applicaiton to triangular billiards}, 
      Invent. Math. {\bf 97} (1990), 117-171.  

\end{thebibliography}
\end{document}